\documentclass[12pt,fleqn]{article}

\usepackage{amsmath}
\usepackage{amssymb}
\usepackage{natbib}
\usepackage{graphicx}
\usepackage{color}
\usepackage{float} 	

\title{ANN-MoC Method for Inverse Transient Transport Problems in One-Dimensional Geometry}

\author
{\small\rm \begin{tabular}{l}
       \textbf{Nelson Garcia Roman} - \texttt{ngroman1992@gmail.com}\\%
       \textbf{Pedro Costa dos Santos} - \texttt{pedro.costa4137@gmail.com}\\
       \textbf{Pedro Henrique de Almeida Konzen} - \texttt{pedro.konzen@ufrgs.br}\\
       {\fontsize{11}{0}\selectfont UFRGS - IME, Porto Alegre, Brazil}
  \end{tabular}}

\date{}

\usepackage{booktabs}
\newcommand{\e}{\mathrm{e}}

\begin{document}

\maketitle

\begin{abstract}
  The inverse problems of particle neutral transport models have many important engineering and medical applications. Safety protocols, quality control procedures, and optical medical solutions can be developed based on inverse transport solutions. In this work, we propose the ANN-MoC method to solve the inverse transient transport problem of estimating the absorption coefficient from measurements of the scalar flux at the boundaries of the model domain. The main idea is to train an Artificial Neural Network (ANN) from data generated by direct solutions computed by a Method of Characteristics (MoC) solver. The direct solver is tested on a problem with a manufactured solution. And, the proposed ANN-MoC method is tested on two inverse problems. In the first, the medium is homogeneous and has a constant absorption coefficient. In the second, a heterogeneous medium is considered, with the absorption coefficient constant by parts. Very accurate ANN estimations have been achieved for these two problems, indicating that the quality of the results relies on the accuracy of the direct solver solutions. The results show the potential of the proposed approach to be applied to more realistic inverse transport problems.

\end{abstract}

\newcommand{\p}{\partial}

\section{INTRODUCTION}
Particle neutral transport phenomena appear in many important engineering and medical applications. The fundamental modeling of radiative heat transfer \citep{Modest2013a} and neutron transport \citep{Lewis1984a} have the linear Boltzmann equation in common. Not restricted to, the first is important in engineering applications at high temperatures, such as glass and ceramic manufactures \citep{Larsen2002a}, combustion chambers \citep{Frank2004a} and solar energy production \citep{Fuqiang2017a}. The other has applications in nuclear energy production \citep{Stacey2007a} and optical medicine \citep{Wang2007a}. In these areas, inverse transport solutions can enhance the development of safety protocols, quality control procedures and, in medicine, contribute to the development of optical medical solutions \citep{Bal2009a}.

We consider the time-dependent linear Boltzmann equation with initial and boundary conditions and with isotropic scattering
\begin{subequations}\label{eq:tp}
  \begin{align}
    &\forall \mu: ~\frac{1}{c}\frac{\p}{\p t}I(t,x,\mu) + \mu\frac{\p I}{\p x} + \sigma_tI = \sigma_s\Psi(t,x) \nonumber\\
    &\qquad\qquad\qquad\qquad\qquad\qquad\qquad+ q(t,x,\mu), ~(t,x)\in (0,t_f]\times\mathcal{D},\label{eq:te}\\
    &\forall \mu:~I(0,x,\mu) = I_0(x,\mu), ~x\in\mathcal{D},\label{eq:tic}\\ 
    &\forall \mu>0:~I(t,a,\mu) = I_{\text{in,a}}(t,\mu), ~t\in (0, t_f],\label{eq:tbca}\\    
    &\forall \mu<0:~I(t,b,\mu) = I_{\text{in,b}}(t,\mu), ~t\in (0, t_f].\label{eq:tbcb}
  \end{align}
\end{subequations}
By considering a radiative transport problem, $I(t,x,\mu)$ $(\text{W}/\text{sr})$ denotes the radiation intensity at time $0 \leq t \leq t_f$ ($\text{s}$) at point $x\in D=[a, b]$ ($\text{m}$), and in the direction $-1 < \mu < 1$, $\mu\neq 0$. The average speed of light in the medium is denoted by $c$ ($\text{m}/\text{s}$). The total absorption coefficient is denoted by $\sigma_t = \kappa + \sigma_s$ $(1/\text{m})$, where $\kappa$ $(1/\text{m})$ and $\sigma_s$ $(1/\text{m})$ are, respectively, the absorption and scattering coefficients. The sources are denoted by $q(t,x,\mu)$ $(\text{W}/(\text{m}\cdot\text{sr}))$ in the domain and $I_{\text{in,a}}(t,\mu)$, $I_{\text{in,a}}(t,\mu)$ $(\text{W}/\text{sr})$ at boundaries. At $t = 0$, initial condition $I = I_0(x,\mu)$ $(\text{W}/\text{sr})$ is assumed. The average scalar flux $(\text{W}/\text{sr})$ is defined as
\begin{equation}\label{eq:sf}
  \Psi(t,x) := \frac{1}{2}\int_{-1}^1I(t,x,\mu)\,d\mu.
\end{equation}

Inverse transport problems have been the subject of important research for many decades. The books of \cite{Ozisik2021a} and \cite{MouraNeto2021a} discuss the fundamental methods applied to the solution of inverse problems. Concerning the problems of parameter estimation, the main approaches consist of estimating parameters as solutions to an associated minimization problem. The problem can then be solved by optimization methods, which usually require a good initial approximation of the solution. When this is not known, meta-heuristic algorithms can be applied to this end (see, for instance, the paper of \cite{Lobato2010a}). Alternatively, Deep Learning \citep{Goodfellow2016a} techniques are also applied \citep{Bokar1999a, Lugon2009a}. A well-known approach is to train an Artificial Neural Network (ANN) \citep{Haykin2008a} with data samples built from solutions to the associated direct problem.

In this context, we propose the ANN-MoC method to solve the inverse transport problem of estimating the absorption coefficient from measurements of the scalar flux at the boundaries of the model domains. The main idea is to train an ANN from data generated by the direct solutions of Eq.~\eqref{eq:tp} computed by a solver based on the Method of Characteristics (MoC) \citep{Evans1998a}. 

In the following, the methodology of the MoC direct solver and the ANN model are presented. Numerical experiments with the proposed approach are then discussed. A problem with a manufactured solution is used as a test case for the direct solver. Then, the ANN-MoC method is applied to solve two inverse problems. In the first, the medium is assumed to have a homogeneous absorption coefficient. In the second, the approach is applied to a heterogeneous test case. Conclusions are then presented.

\section{THE ANN-MOC METHOD}
The ANN-MoC method consists of solving the inverse transport problem by an Artificial Neural Network  (ANN) trained from data generated by the direct solving of a set of transport problems by the Method of Characteristics (MoC).

\subsection{MoC direct solver}
The direct solver computes an approximation of Eq.~\eqref{eq:tp} built with the Discrete Ordinates Method (DOM) \citep{Modest2013a} followed by an implicit Euler time discretization \citep{Stoer2000a}. The raised system of ordinary differential equations is decoupled by a Source Iteration (SI) \citep{Modest2013a} scheme and, then, solved with the Method of Characteristics (MoC) \citep{Evans1998a}.

\vspace{0.5cm} 

\textbf{\textit{Discrete ordinates formulation.}} The following DOM form of Eq.~\eqref{eq:tp} is obtained by assuming the Gauss-Legendre quadrature $\left\{(\mu_j,\omega_j)\right\}_{j=1}^{n_q}$, with even $n_q>1$,
\begin{subequations}\label{eq:tp_dom}
  \begin{align}
    &\forall \mu_j: ~\frac{1}{c}\frac{\p}{\p t}I_j(t,x) + \mu_j\frac{\p I_j}{\p x} + \sigma_tI_j = \sigma_s\Psi(t,x) \nonumber\\
    &\qquad\qquad\qquad\qquad\qquad\qquad\qquad+ q_j(t,x), ~(t,x)\in (0,t_f]\times\mathcal{D},\label{eq:te_dom}\\
    &\forall \mu_j:~I_j(0,x) = I_{j,0}(x), ~x\in\mathcal{D},\label{eq:tic_dom}\\ 
    &\forall \mu_j>0:~I_j(t,a) = I_{\text{j,in,a}}(t), \forall t\in (0, t_f],\label{eq:tbca_dom}\\    
    &\forall \mu_j<0:~I_j(t,b) = I_{\text{j,in,b}}(t), \forall t\in (0, t_f],\label{eq:tbcb}
  \end{align}
\end{subequations}
where the notation $I_j(t,x) \approx I(t,x,\mu_j)$ (analogous to the others) is assumed with $j = 1, 2, \dotsc, n_q$. The scalar flux is approximated by
\begin{equation}
  \Psi(t,x) \approx \frac{1}{2}\sum_{j=1}^{n_q} I_j\omega_j.
\end{equation}

\vspace{0.5cm} 

\textbf{\textit{Time discretization.}} For the time discretization, it is assumed that $t^{(k)} = kh_t$, $k = 0, 1, 2, \dotsc, n_t$, $h_t = t_f/n_t$. The implicit Euler formulation of Eq.~\eqref{eq:tp_dom} gives an iterative procedure with initialization
\begin{equation}
  \forall \mu_j: I_j^{(0)}(x) = I_{j,0}(x), ~x\in\mathcal{D}, \label{eq:tic_euler}
\end{equation}
and the following steps
\begin{subequations}\label{eq:tp_euler}
  \begin{align}
    &\forall \mu_j: ~\frac{1}{c}\frac{I_j^{(k+1)} - I_j^{(k)}}{h_t} + \mu_j\frac{\p I_j^{(k+1)}}{\p x} + \sigma_tI_j^{(k+1)}(x) = \sigma_s\Psi^{(k+1)}(x) \nonumber\\
    &\qquad\qquad\qquad\qquad\qquad\qquad\qquad\qquad\qquad\qquad\;+ q_j^{(k+1)}(x), \label{eq:te_euler}\\
    &\forall \mu_j>0:~I_j^{(k+1)}(a) = I^{(k+1)}_{j,\text{in},a}, \label{eq:tbca_euler}\\    
    &\forall \mu_j<0:~I_j^{(k+1)}(b) = I^{(k+1)}_{j,\text{in},b}, \label{eq:tbcb_euler}
  \end{align}
\end{subequations}
where the notation $I_j^{(k)}(x) \approx I\left(t^{(k)},x,\mu_j\right)$ (analogous to the others) is assumed with $k = 0, 1, 2, \dotsc, n_t-1$. For the sake of simplicity, in the following the index $k$ will be suppressed, with $I^{(1)}_j$ denoting $I^{(k+1)}_j$ and $I^{(0)}_j = I^{(k)}_j$ (analogous to the others).

\vspace{0.5cm} 

\textbf{\textit{Source iteration.}} The decoupling of system Eq.~\eqref{eq:tp_euler} is performed with the Source Iteration (SI) technique. From a given initial scalar flux $\Psi^{(0,0)}$, successive approximations $\Psi^{(1,l)}$ are iteratively computed from
\begin{subequations}\label{eq:tp_si}
  \begin{align}
    &\forall \mu_j: ~\frac{1}{c}\frac{I_j^{(1,l+1)} - I_j^{(0)}}{h_t} + \mu_j\frac{\p I_j^{(1,l+1)}}{\p x} + \sigma_tI_j^{(1,l+1)}(x) = \sigma_s\Psi^{(1,l)}(x) \nonumber\\
    &\qquad\qquad\qquad\qquad\qquad\qquad\qquad\qquad\qquad\qquad\quad\;+ q_j^{(1)}(x), \label{eq:te_si}\\
    &\forall \mu_j>0:~I_j^{(1,l+1)}(a) = I^{(1,l+1)}_{j,\text{in},a}, \label{eq:tbca_si}\\    
    &\forall \mu_j<0:~I_j^{(1,l+1)}(b) = I^{(1,l+1)}_{j,\text{in},b}, \label{eq:tbcb_si}
  \end{align}
\end{subequations}
where
\begin{equation}
  \Psi^{(1,l)} := \frac{1}{2}\sum_{j=1}^{n_q}I^{(1,l)}_j\omega_j,
\end{equation}
for $l = 0, 1, 2, \ldots$ until some given stop criteria is fulfilled.

\vspace{0.5cm} 

\textbf{\textit{Method of characteristics.}} At each time step and each source iteration, Eq.~\eqref{eq:tp_si} is solved by the Method of Characteristics (MoC). First, it is observed that Eq.~\eqref{eq:te_si} can be rewritten as
\begin{equation}\label{eq:te_aux}
  \mu_j\frac{\p I_j^{(1,l+1)}}{\p x} + \left(\sigma_t + \frac{1}{c h_t}\right)I_j^{(1,l+1)}(x) = \sigma_s\Psi^{(1,l)}(x) + q_j^{(1)}(x) + \frac{1}{c h_t}I_j^{(0)}
\end{equation}
Again, for the sake of simplicity, the index $j$ is suppressed in the following.

The MoC form of Eq.~\eqref{eq:te_aux} is obtained by assuming $x(s) = x_0 + s\mu$, $s\in\mathbb{R}$, from where Eq.~\eqref{eq:te_aux} is rewritten as
\begin{equation}\label{eq:te_moc}
  \frac{d}{d s}I^{(1,l+1)} + \left(\sigma_t + \frac{1}{c h_t}\right)I^{(1,l+1)}(s) = \sigma_s\Psi^{(1,l)}(s) + q_j^{(1)}(s) + \frac{1}{c h_t}I_j^{(0)}(s).
\end{equation}
This linear first-order differential equation can now be solved using an integration factor, which gives the solution form
\begin{equation}\label{eq:moc_sol}
  I^{(1,l+1)}(s) = I^{(1,l+1)}(0)e^{-\int_0^{s}\tilde{\sigma}_t\,ds'} + \int_0^s S^{(l)}(s')e^{-\int_{s'}^{s}\tilde{\sigma}_t\,ds''}\,ds',
\end{equation}
where $\tilde{\sigma}_t := \sigma_t + 1/(c h_t)$ and
\begin{equation}
  S^{(l)}(s) := \sigma_s\Psi^{(1,l)}(s) + q_j^{(1)}(s) + \frac{1}{c h_t}I_j^{(0)}(s).
\end{equation}
One observes that, choosing $x_0 = a$, Eq.~\eqref{eq:moc_sol} gives the particle intensity $I^{(1,l+1)}(x(s))$ at each domain point $x(s)$ for a given direction $\mu > 0$. Analogously, choosing $x_0=b$, one obtains the particle intensity at any domain point for a given direction $\mu < 0$.

\vspace{0.5cm} 

\textbf{\textit{Direct solver algorithm.}} Assuming a spatial mesh with $n_x$ nodes $x_i = a + ih_x$, and mesh size $h_x = (b-a)/n_x$, $i = 0, 1, 2, \dotsc, n_x$, the direct solver algorithm can be summarized as follows:
\begin{itemize}
\item[1.] Set time, mesh and quadrature parameters.
\item[2.] From initial condition, set
  \begin{subequations}
    \begin{align}
      &I_{i,j}^{(0)} \leftarrow I(0,x_i,\mu_j), ~\forall i,j,\\
      &\Psi^{(0)}_i \leftarrow \frac{1}{2}\sum_{j=1}^{n_q}I_{i,j}^{(0)}\omega_j, ~\forall i.
    \end{align}
  \end{subequations}
\item[3.] (Time loop.) For $k = 0, 1, 2, \dotsc, n_t$:
  \begin{itemize}
  \item[a.] (SI loop.) For $l = 0, 1, 2, \dotsc, n_{\text{s.i.}}$:
    \begin{itemize}
    \item[a.1.] For $j = 1, 2, \dotsc, n_q$ and $\mu_j>0$:
      \begin{itemize}
      \item For $i = 0, 1, 2, \dotsc, n_x-1$:
        \begin{equation}
          I_{i+1,j}^{(1,l+1)} \leftarrow I_{i,j}^{(1,l+1)}e^{-\int_0^{s}\tilde{\sigma}_t\,ds'} + \int_0^s S^{(l)}(s')e^{-\int_{s'}^{s}\tilde{\sigma}_t\,ds''}\,ds'
        \end{equation}
      \end{itemize}
    \item[a.2.] For $j = 1, 2, \dotsc, n_q$ and $\mu_j<0$:
      \begin{itemize}
      \item For $i = n_x, n_x-1, \dotsc, 1$:
        \begin{equation}
          I_{i-1,j}^{(1,l+1)} \leftarrow I_{i,j}^{(1,l+1)}e^{-\int_0^{s}\tilde{\sigma}_t\,ds'} + \int_0^s S^{(l)}(s')e^{-\int_{s'}^{s}\tilde{\sigma}_t\,ds''}\,ds'
        \end{equation}
      \end{itemize}
    \item[a.3.] Compute new scalar flux
      \begin{equation}
        \Psi^{(l+1)}_i \leftarrow \frac{1}{2}\sum_{j=1}^{n_q}I_{i,j}^{(1,l+1)}\omega_j, ~\forall i.
      \end{equation}
    \item[a.4.] SI stop criterion.
    \end{itemize}
  \end{itemize}
\end{itemize}

\subsection{ANN inverse model}
The inverse problem is solved by training a Multilayer Perceptron neural network (MLP) \citep{Haykin2008a} from a data set $\{(\pmb{\psi}^{(s)}_{\text{train}}, \kappa^{(s)}_{\text{train}})\}_{s=1}^{n_{\text{train}}}$ computed by solving the direct problem for several values of the absorption coefficient. The MLP of $n_l$ layers can be written as
\begin{equation}
  \tilde{\kappa} = \mathcal{N}\left(\pmb{\psi}; \left\{\left(\pmb{f}^{(l)},\pmb{b}^{(l)}, W^{(l)}\right)\right\}_{l=1}^{n_l}\right),
\end{equation}
where, at the $l$-th network layer with $m_l$ neuron units, $\left(\pmb{f}^{(l)},\pmb{b}^{(l)}, W^{(l)}\right)$ denotes the triple of the activation function, the bias $m_l$-vector, and the weights $m_l\times m_{l-1}$-matrix. By denoting the input $\pmb{y}^{(0)} = \pmb{\psi}$ of detector measurements, its forward propagation through the network layers $l = 1, 2, \dotsc, n_l$ is given by
\begin{equation}
  \pmb{y}^{(l+1)} = \pmb{f}^{(l+1)}\left(W^{(l+1)}\pmb{y}^{(l)} + \pmb{b}^{(l+1)}\right),
\end{equation}
and the output is the estimated absorption coefficient $\tilde{\kappa} = y^{(n_l)}$.

\vspace{0.5cm} 

\textbf{\textit{Basic training algorithm.}} The basic training algorithm can be summarized as follows:
\begin{itemize}
\item[1.] Set the MLP architecture.

  Sets the number of input units $m_0$, the number of units at each hidden layer $m_l$, its activation functions $\pmb{f}^{(l)}$, and initial biases $\pmb{b}^{(l)}$, weights $W^{(l)}$ and a global learning rate $l_r>0$.
\item[2.] Loop over epochs $e \leftarrow 1, 2, \dotsc, n_{\text{epochs}}$:
  \begin{itemize}
  \item[2.a.] Forward the training set.
    \begin{equation}
      \tilde{\pmb{\kappa}}_{\text{train}} \leftarrow \mathcal{N}\left(\pmb{\psi}_{\text{train}}\right)
    \end{equation}
  \item[2.b.] Compute the loss function.
    \begin{equation}
      \mathcal{L} \leftarrow \frac{1}{n_{\text{train}}}\sum_{s=1}^{n_{\text{train}}}\left|\tilde{\kappa}_{\text{train}}^{(s)} - \kappa_{\text{train}}^{(s)}\right|^2
    \end{equation}
  \item[2.c.] Backward the loss function to compute the gradients.
    \begin{equation}
      \frac{\p\mathcal{L}}{\p\pmb{b}^{(l)}}, ~\frac{\p\mathcal{L}}{\p W^{(l)}}, ~l = 1, 2, \dotsc, n_l.
    \end{equation}
  \item[2.d.] Perform an optimizer gradient based step.
    \begin{subequations}
      \begin{align}
        &W^{(l)} \leftarrow W^{(l)} - l_{r}\frac{\p\mathcal{L}}{\p W^{(l)}},\\
        &\pmb{b}^{(l)} \leftarrow \pmb{b}^{(l)} - l_{r}\frac{\p\mathcal{L}}{\p \pmb{b}^{(l)}}
      \end{align}
    \end{subequations}
  \item[2.e.] Check the stopping criterion.
  \end{itemize}
\end{itemize}

\vspace{0.5cm} 

\textbf{\textit{ANN model test}}. The test of the trained neural network model consists of verifying its performance for a new data set $\{(\pmb{\psi}^{(s)}_{\text{test}}, \kappa^{(s)}_{\text{test}})\}_{s=1}^{n_{\text{test}}}$ which has not been used for training. The test data set has also been computed by solving the direct problem for several values of the absorption coefficient. The accuracy of the network estimated values $\tilde{\kappa}^{(s)}_{\text{test}}$ can be measured by the mean squared error $\mathcal{L}_{\text{test}}$ and the coefficient of determination.

\section{RESULTS}
Tests of the proposed ANN-MoC approach are now presented. First, the MoC direct solver is tested on a transport problem with a manufactured solution. Two inverse problems are then discussed. In the first, a homogeneous medium is assumed, and in the second a two region heterogeneous medium is considered.

\subsection{Direct solver test}
In order to test the direct solver, we have considered the manufactured solution
\begin{equation}
  \hat{I}(t,x,\mu) := e^{-\sigma_t|x-t|^2}.
\end{equation}
By substituting $\hat{I}$ into Eq.~\eqref{eq:te}, the source is found to be
\begin{equation}
  q(t,x,\mu) = \left[2\sigma_t(1-\mu)(x-t)+\kappa\right]e^{-\sigma_t|x-t|^2},
\end{equation}
and, from the definition of the scalar flux Eq.~\eqref{eq:sf}, one has $\hat{\Psi} = \hat{I}$.

After numerical tests, we have chosen the solver parameters $h_t = 0.01$, $n_x = n_q = 100$, and $\texttt{tol}=1.49\times 10^{-8}$ as the absolute $L^2$-norm tolerance for the SI stopping criterion. Table~\ref{tab:ds_test} shows a comparison between the direct solver approximations and the exact scalar flux solutions at $t_f=1.0$ for different absorption coefficients. The relative $L^2$-error is denoted by $\varepsilon_{\text{rel}}$ and indicates that the chosen parameters were enough for the direct solver to produce an accurate solution with $\varepsilon_{\text{rel}} < 10^{-2}$.

\begin{table}[H]
  \centering
  \caption{Comparison between the direct solver approximations and the exact solution at $t_f=1.0$.}
  \vspace{12pt}
  \begin{tabular}{r|ccc|c}\toprule
    $\kappa$ & $\Psi(0.0)$ & $\Psi(0.5)$ & $\Psi(1.0)$ & $\varepsilon_{\text{rel}}$\\\midrule
    $0.9$    & $3.667\e-1$ & $7.748\e-1$ & $9.974\e-1$ & $4.5\e-3$\\
    $0.5$    & $3.664\e-1$ & $7.740\e-1$ & $9.971\e-1$ & $5.3\e-3$\\
    $0.1$    & $3.660\e-1$ & $7.730\e-1$ & $9.968\e-1$ & $6.4\e-3$\\\midrule
    exact    & $3.679\e-1$ & $7.788\e-1$ & $1.000\e+0$\\\bottomrule
  \end{tabular}
  \label{tab:ds_test}
\end{table}

\subsection{Inverse problem 1 -  homogeneous medium}
As a first test case, we consider a homogeneous medium with a constant absorption coefficient. The inverse problem consists of estimating $0.1 < \kappa < 0.9$ from detectors measurements of the scalar fluxes on the boundaries of the domain $\mathcal{D} = [0, 1]$ and at time $t_{d,3}=3.0$. Boundary conditions are taken as $I(t,0,\mu) = 1$, for all $\mu >0$, and $I(t,1,\mu) = 0$, for all $\mu<0$. The null source is considered, and the initial condition is $I(0,0,\mu) = 1$, $\mu>0$, and $I(0,x,\mu) = 0$ for all $x>0$.

\begin{figure}[H]
  \centering
  \includegraphics[width=0.39\textwidth]{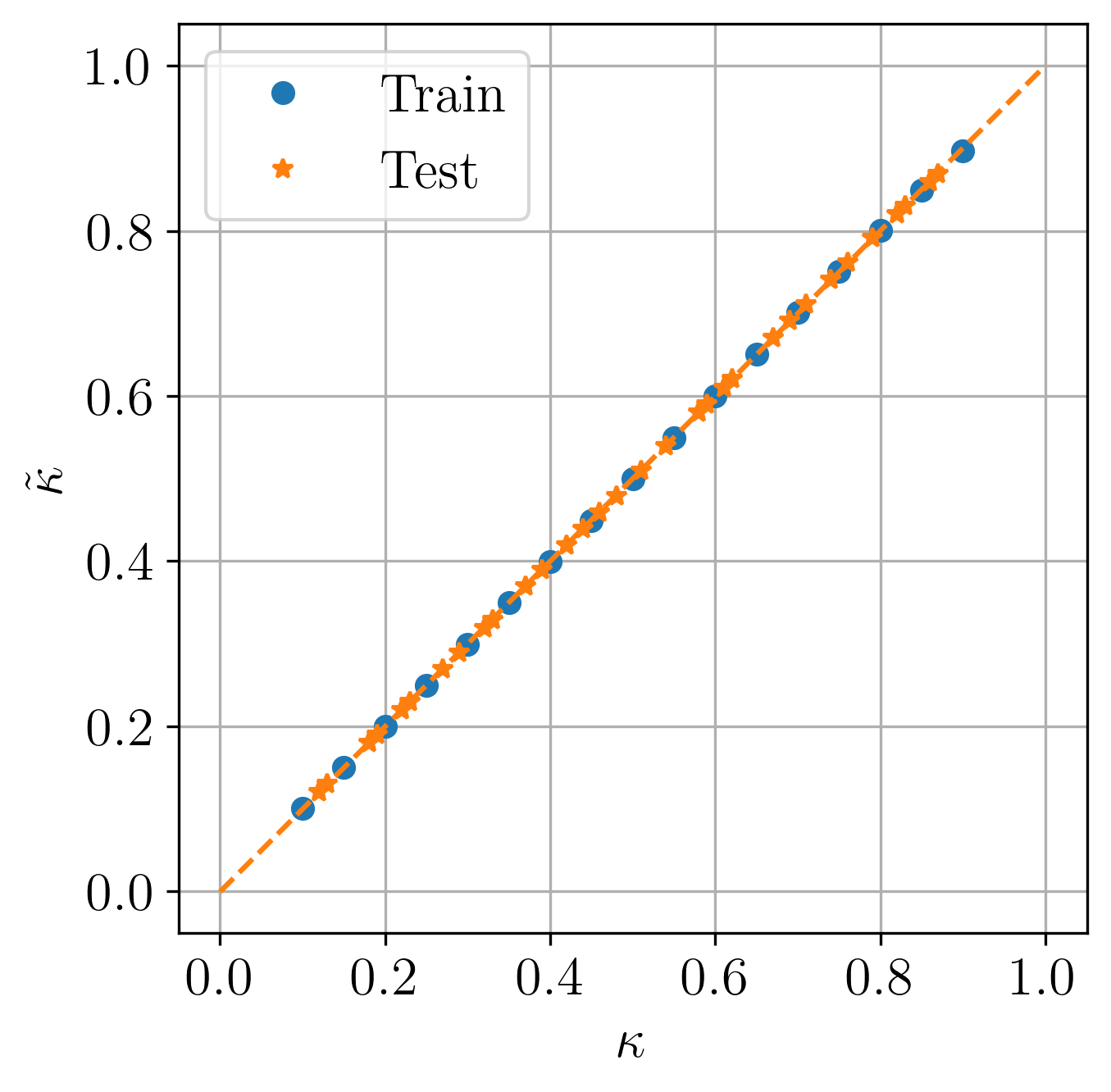}~
  \includegraphics[width=0.39\textwidth]{./data/prob2_res/fig}
  \caption{Inverse problem 1. Left: training (circles) and test (stars) samples. Right: expected $\kappa$ \textit{versus} estimated $\tilde{\kappa}$.}
  \label{fig:prob2}
\end{figure}

The ANN inverse model has the detectors measurements $d_0 = \Psi(t_{d,3},0)$, $d_1 = \Psi(t_{d,3},1)$ as inputs and outputs the estimated absorption coefficient $\tilde{\kappa}$. For its training, we have used the direct solver to build a training set $\{(\pmb{d}^{(s)}_{\text{train}}, \kappa^{(s)}_{\text{train}})\}_{s=1}^{n_{\text{train}}}$ of $n_{\text{train}} = 17$ samples (patterns) with $\kappa^{(s)} = 0.1 + (s-1)h_s$, $h_s=0.05$. The test set $\{(\pmb{d}^{(s)}_{\text{test}}, \kappa^{(s)}_{\text{test}})\}_{s=1}^{n_{\text{test}}}$ has been generated with $n_{\text{test}} = 32$ with random choices $0.1 < \kappa^{(s)} < 0.9$ (see Fig.~\ref{fig:prob2}, left).

By following a trial and error strategy, we choose a MLP model with architecture $2-25-25-25-1$ (tow inputs, three hidden layers with 25 neurons each, and one output neuron), the hyperbolic tangent and the identity as activation functions in the hidden and output layers, respectively. With about $n_{\text{epoch}}=6000$ the model reaches a mean squared error $\mathcal{L}_{\text{train}} < 10^{-6}$ and the coefficient of determination $R^2_{\text{train}} = 0.99998$. The application of the trained model to the test data gave results with $\mathcal{L}_{\text{test}} < 10^{-6}$ and $R^2_{\text{test}} = 0.99999$ (see Fig.~\ref{fig:prob2}, right).

The relatively high number of epochs required for training convergence can be explained by the lack of data pre-processing and an early-stopping method. However, subject of further developments, this was not a current concern since the training ended in a few seconds of computation and no over-fitting (or over-training) was observed.

\subsection{Inverse problem 2 - heterogeneous medium}
In the second test case, we consider a heterogeneous medium with piece-wise constant absorption coefficients
\begin{equation}
  \kappa(x) = \left\{
    \begin{array}{ll}
      \kappa_1 &, 0 \leq x \leq 0.5,\\
      \kappa_2 &, 0.5 < x \leq 1.
    \end{array}
\right.
\end{equation}

The inverse problem consists in estimating $0.1 < \kappa_1,\kappa_2 < 0.9$ from detectors measurements of the scalar fluxes on the boundaries of the domain $\mathcal{D} = [0, 1]$ and at the times $t_{d,2}=2.0$ and $t_{d,3}=3.0$. The initial and boundary conditions, as well as the source are the same as for inverse problem 1.

The ANN inverse model has the detectors measurements $\pmb{d}_0 = \left(\Psi(t_{d,2},0), \Psi(t_{d,3},0)\right)$, $d_1 = \left(\Psi(t_{d,2},1), \Psi(t_{d,3},1)\right)$ as inputs and outputs the estimated absorption coefficients $\tilde{\kappa}_1$ and $\tilde{\kappa}_2$. For its training, we have used the direct solver to compute the training set $\{(\pmb{d}^{(s)}_{\text{train}}, \pmb{\kappa}^{(s)}_{\text{train}})\}_{s=1}^{n_{\text{train}}}$ of $n_{\text{train}} = 81$ samples (patterns) with $\kappa_{1,2}^{(s)} = 0.1 + (s-1)h_s$, $h_s=0.1$. The test set $\{(\pmb{d}^{(s)}_{\text{test}}, \pmb{\kappa}^{(s)}_{\text{test}})\}_{s=1}^{n_{\text{test}}}$ has been generated with $n_{\text{test}} = 64$ random choices $0.1 < \kappa_{1,2}^{(s)} < 0.9$.

\begin{figure}[H]
  \centering
  \includegraphics[width=0.49\textwidth]{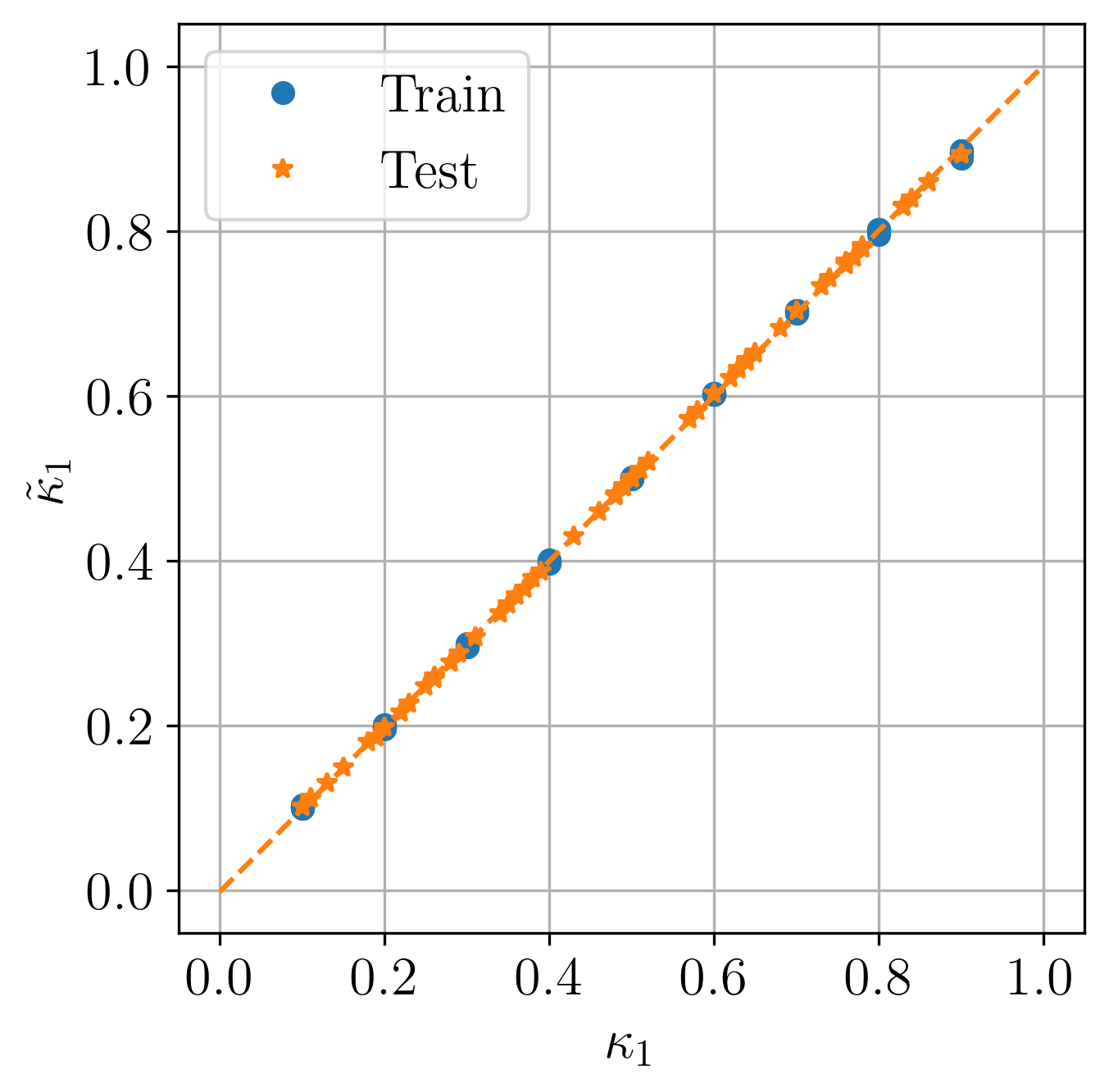}~
  \includegraphics[width=0.49\textwidth]{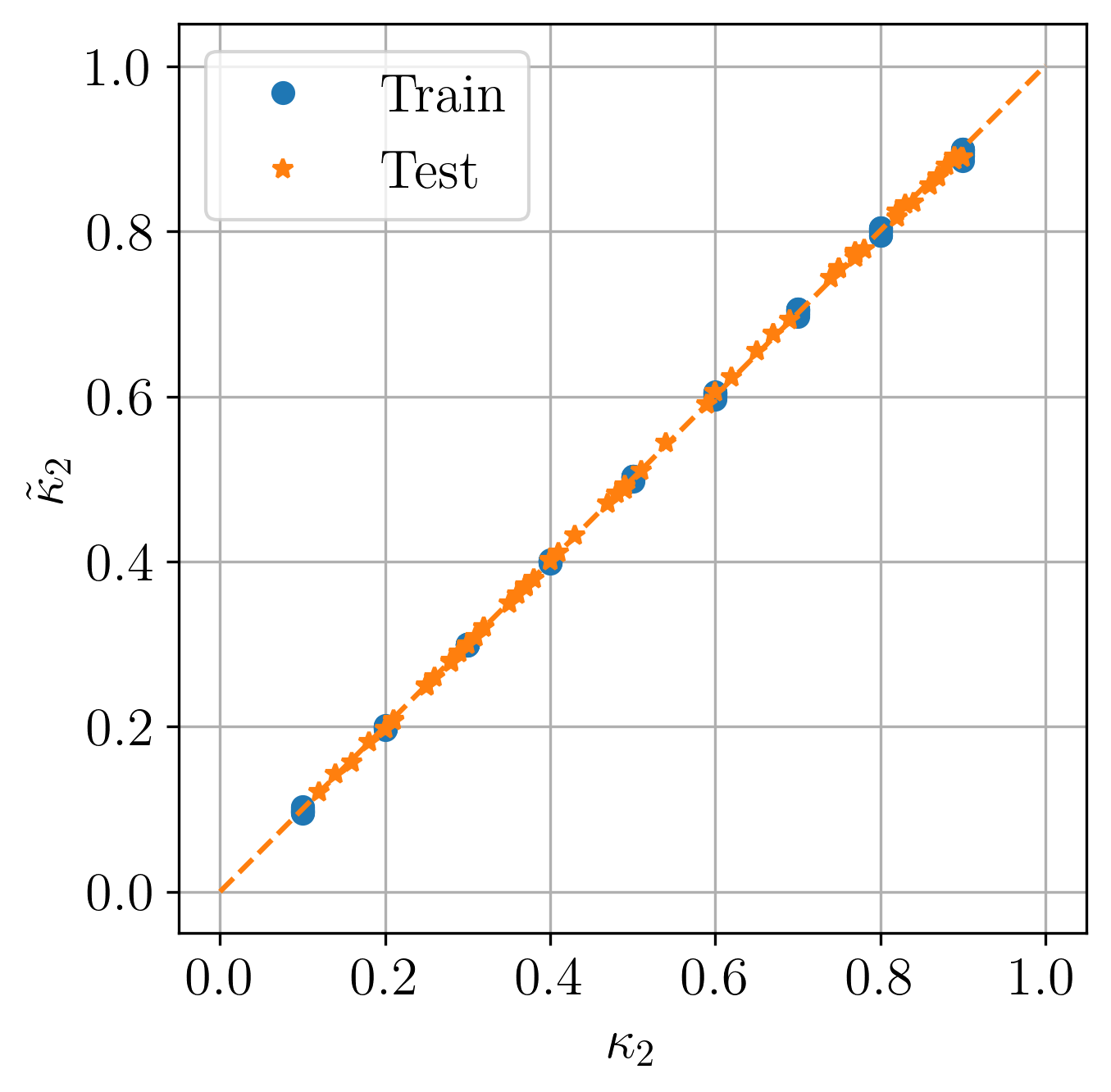}
  \caption{Inverse problem 2: expected $\pmb{\kappa}$ \textit{versus} estimated $\tilde{\pmb{\kappa}}$. Left: $\kappa_1$. Right: $\kappa_2$.}
  \label{fig:prob3}
\end{figure}

By following a trial and error strategy, we choose an ANN model with architecture $4-25-25-25-25-2$, the hyperbolic tangent and the identity as activation functions in the hidden and output layers, respectively. With about $n_{\text{epoch}}=2800$ the model reaches a mean squared error $\mathcal{L}_{\text{train}} < 10^{-5}$ and the coefficient of determination $R^2_{\text{train}} = 0.9999$. The application of the trained model to the test data gave results with $\mathcal{L}_{\text{test}} < 10^{-5}$ and $R^2_{\text{test}} = 0.9999$. Figure~\ref{fig:prob3} shows the expected \textit{versus} estimated absorption coefficients for the training and test samples.

\section{CONCLUSIONS}
In this work, we present the ANN-MoC method to solve the inverse transport problem of estimating the absorption coefficient from measurements of the scalar fluxes at the boundaries of the computational domain. The main idea is to train an Artificial Neural Network (ANN) from data computed by direct solutions of a set of transport problems based on the Method of Characteristics (MoC).

The reported applications of two different inverse transport problems show the potential of the proposed approach. The ANN models could very efficiently learn from the generated data and estimate the absorption coefficients with high accuracy. This indicates that the accuracy of the inverse model highly depends on the accuracy of the direct solutions. Further developments should aim to improve the direct solver. Improvements in the accuracy of the solution and, mainly, in the computational performance are important to feed the ANN model with a better-quality data set. Solutions to more complex inverse transport problems could also profit from the proposed approach, but again, it will demand further improvements of the direct solver. Finally, the use of the proposed methodology for realistic problems depends on how good the direct transport model is for the intended application.

\subsection*{\textit{Acknowledgements}}
This study was financed in part by the Coordenação de Aperfeiçoamento de Pessoal de Nível Superior – Brasil (CAPES) – Finance Code 001.


\end{document}